\documentclass[10pt]{article}
\usepackage{lmodern}
\usepackage{amsmath}
\usepackage{authblk}
\usepackage{amsfonts}
\usepackage{amssymb}
\usepackage[english]{babel}
\usepackage{color}
\usepackage{graphicx}

\usepackage[latin1]{inputenc}
\usepackage{tikz}
\tikzset{liltext/.style={font=\tiny}}
\usetikzlibrary{decorations.pathreplacing,calc}
\usepackage{rotating}
\usepackage{float}

\usepackage{amsthm}
\usepackage{graphicx}
\usepackage{mathrsfs}
\usepackage{makecell}
\usepackage{microtype}
\usepackage{mathscinet}
\usepackage{array}
\usepackage{multirow}
\usepackage{enumerate}
\usepackage{hyperref}
\hypersetup{hidelinks}
\newtheorem{defn}{Definition}[section]
\newtheorem{theorem}{Theorem}[section]
\newtheorem{prop}{Proposition}[section]

\newtheorem{remark}{Remark}[section]

\newcommand{\R}{\mathbb{R}}

\newcommand{\N}{\mathbb{N}}
\newcommand{\e}{\varepsilon}
\newcommand{\ity}{\infty}
\newcommand{\ml}{\mathcal}
\newcommand{\mb}{\mathbb}
\DeclareMathOperator{\lin}{lin}
\DeclareMathOperator{\non}{non}

\title{Sharp lifespan estimates for the weakly coupled system of semilinear damped wave equations in the critical case}
\author[1]{Wenhui Chen}
\affil[1]{School of Mathematical Sciences, Shanghai Jiao Tong University, 200240 Shanghai, China}
\author[2,3]{Tuan Anh Dao\thanks{Corresponding author: Tuan Anh Dao (anh.daotuan@hust.edu.vn)}}
\affil[2]{School of Applied Mathematics and Informatics, Hanoi University of Science and Technology, No.1 Dai Co Viet road, Hanoi, Vietnam}
\affil[3]{Institute of Mathematics, Vietnam Academy of Science and Technology, No.18 Hoang Quoc Viet road, Hanoi, Vietnam}
\date{}
\begin{document}
		\maketitle
\begin{abstract}
The open question, which seems to be also the final part, in terms of studying the Cauchy problem for the weakly coupled system of damped wave equations or reaction-diffusion equations, is so far known as the sharp lifespan estimates in the critical case. In this paper, we mainly investigate lifespan estimates for solutions to the weakly coupled system of semilinear damped wave equations in the critical case. By using a suitable test function method associated with nonlinear differential inequalities, we catch upper bound estimates for the lifespan. Moreover, we establish polynomial-logarithmic type time-weighted Sobolev spaces to obtain lower bound estimates for the lifespan in low spatial dimensions. Then, together with the derived lifespan estimates, new and sharp results on estimates for the lifespan in the critical case are claimed. Finally, we give an application of our results to the semilinear reaction-diffusion system in the critical case.
	\medskip\\
	\textbf{AMS Classification (2020)} Primary 35G55; Secondary 35L71. 
\end{abstract}
\fontsize{12}{15}
\selectfont


\section{Introduction}\label{Sec_Intro}
In the present paper, we are interested in exploring sharp lifespan estimates for the weakly coupled system of semilinear classical damped wave equations in the critical case, namely,
\begin{align}\label{Eq_Coupled_Damped_Waves}
\begin{cases}
u_{tt}-\Delta u+u_t=|v|^p,&x\in\mb{R}^n,\ t\in(0,T),\\
v_{tt}-\Delta v+v_t=|u|^q,&x\in\mb{R}^n,\ t\in(0,T),\\
(u,u_t,v,v_t)(0,x)=(\varepsilon u_0,\varepsilon u_1,\varepsilon v_0,\varepsilon v_1)(x),&x\in\mb{R}^n,
\end{cases}
\end{align}
where the power exponents $p,q>1$ satisfy the following critical condition:
\begin{align}\label{Critical_Condition}
\alpha_{\max}(p,q):=\frac{\max\{p,q\}+1}{pq-1}=\frac{n}{2}
\end{align}
for any $n\geqslant 1$, $T>0$ and the positive constant $\varepsilon$ describes the size of initial data. Under this critical condition, every non-trivial local (in time) weak solution blows up in finite time. For this reason, a natural question is that whether or not one can describe more detailed information of the lifespan. More specifically, our main motivation of this paper is to report sharp estimates for the lifespan $T=T_\varepsilon$ of solutions to the weakly coupled system \eqref{Eq_Coupled_Damped_Waves} under the critical condition \eqref{Critical_Condition}. Here, the lifespan $T_\varepsilon$ of solutions is understood as the quantity defined by
\begin{align}\label{Lifespan_Defn}
T_\varepsilon:= \sup& \left\{T\in (0,\ity) : \mbox{there exists a unique local (in time) solution} (u,v) \ \mbox{to \eqref{Eq_Coupled_Damped_Waves}}\right.\notag\\
	&\left. \qquad\qquad\qquad \mbox{on } [0,T) \mbox{ with a fixed parameter }\varepsilon>0\right\}.	
\end{align}
 There are a lot of related works begun from 1995 in the original paper \cite{Li-Zhou-1995} in terms of the study of the Cauchy problem \eqref{Eq_Coupled_Damped_Waves}. However, to the best knowledge of authors, the sharp lifespan estimates in the critical case \eqref{Critical_Condition} are completely open even for the semilinear weakly coupled reaction-diffusion systems (see \cite{Fuji-Iked-Waka-2020}). We will partially give answers to the above question in some spatial dimensions by the following sharp estimates:
\begin{align}\label{Sharp_Lifespan}
	T_{\varepsilon}\sim\begin{cases}
		\exp\left( C\varepsilon^{-(p-1)} \right)&\mbox{if}\ \ p=q,\\
		\exp\left(C\e^{-(pq-p_{\mathrm{Fuj}}(n))}\right) &\mbox{if}\ \ p\neq q,
	\end{cases}
\end{align}
where $C>0$ is a constant independent of $\varepsilon$. Here, $p_{\mathrm{Fuj}}(n):=1+2/n$ stands for the well-known Fujita exponent.

Let us now recall several historical background related to our model \eqref{Eq_Coupled_Damped_Waves}. Over the recent decades, the Cauchy problem for semilinear damped wave equation
\begin{align}\label{Eq_Single_Semilinear_Damped_Wave}
	\begin{cases}
		u_{tt}-\Delta u+u_t=|u|^p,&x\in\mb{R}^n,\ t\in(0,T),\\
		(u,u_t)(0,x)=(\varepsilon u_0,\varepsilon u_1)(x),&x\in\mb{R}^n,
	\end{cases}
\end{align}
with $p>1$ has been widely studied. The critical exponent, which is the threshold between global (in time) existence of small data weak solution and blow-up of solutions even for small data, to the single semilinear damped wave equation \eqref{Eq_Single_Semilinear_Damped_Wave} is the so-called \emph{Fujita exponent}
\begin{align*}
	p_{\mathrm{Fuj}}(n):=1+\frac{2}{n}
\end{align*}
for any $n\geqslant 1$. We should remember that the Fujita exponent is well-known as the critical exponent for the corresponding semilinear heat equation (see \cite{Fujita-1966} and references therein)
\begin{align*}
	\begin{cases}
		w_t-\Delta w=|w|^p,&x\in\mb{R}^n,\ t\in(0,T),\\
		w(0,x)= \varepsilon w_0(x),&x\in\mb{R}^n.
	\end{cases}
\end{align*}
Motivated by diffusion phenomenon (see, for example, \cite{Nish-2013,Ikehata-Nishi-2003}), which is a bridge between the asymptotic behavior of solutions to the damped wave equation and that of solutions to the heat equation, one proved that these corresponding critical exponents coincide. In particular, concerning the critical exponent for the semilinear damped wave equation \eqref{Eq_Single_Semilinear_Damped_Wave}, we refer the interested readers to \cite{Zhan-2001,Todo-Yord-2001,Nishihara-2003,Nara-2004,Ikeh-Tani-2005} and references therein. Additionally, to derive the critical regularity of nonlinear terms for the semilinear damped wave equation, the authors in \cite{Eber-Gira-Reis-2020} considered the equation of \eqref{Eq_Single_Semilinear_Damped_Wave} with nonlinearities $\omega(|u|)|u|^{p_{\mathrm{Fuj}}(n)}$ on the right-hand side, where $\omega$ stands for a suitable modulus of continuity. According to the works \cite{Kira-Qafs-2002,Nish-2011,Iked-Ogaw-2016,Lai-Zhou-2019,Fuji-Iked-Waka-2019,Iked-Waka-2020}, the sharp lifespan estimates for all spatial dimensions have been investigated, particularly,
\begin{align*}
	T_{\varepsilon}\sim \begin{cases}
		C\varepsilon^{-\frac{2(p-1)}{2-n(p-1)}}&\mbox{if}\ \ 1<p<p_{\mathrm{Fuj}}(n),\\
		\exp\left(C\varepsilon^{-(p-1)}\right)&\mbox{if}\ \ p=p_{\mathrm{Fuj}}(n),
	\end{cases}
\end{align*}
where $C=C(n,p,u_0,u_1)$ is a positive constant independent of $\varepsilon$. In some sense, the story of the sharp lifespan estimates for the single semilinear damped equation \eqref{Eq_Single_Semilinear_Damped_Wave} has been completed in 2019 by the paper \cite{Lai-Zhou-2019}.

Let us turn to the weakly coupled system of semilinear damped wave equations \eqref{Eq_Coupled_Damped_Waves}. The critical condition \eqref{Critical_Condition}, the so-called critical curve in the $p-q$ plane, to our model \eqref{Eq_Coupled_Damped_Waves} can be described by the following figure:
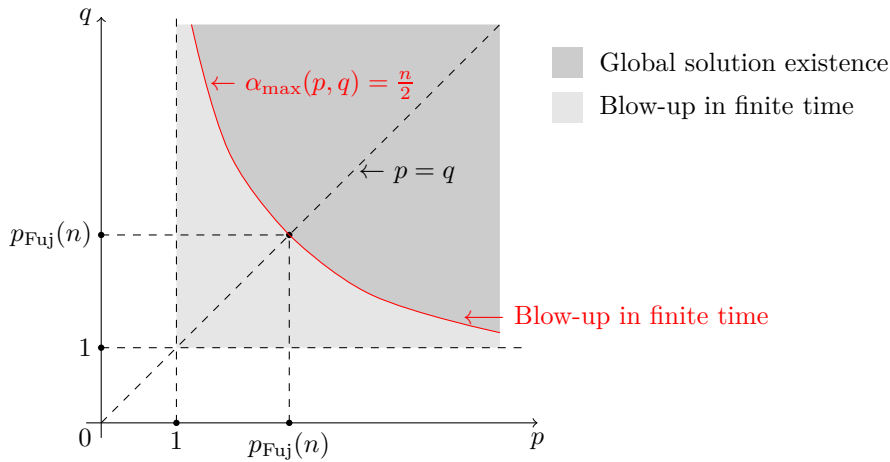
\begin{figure}[http]
	\centering
	\begin{tikzpicture}
	\fill[domain=2:4.8,color=black!10!white] plot[smooth, tension=.7] coordinates {(1.2,5.3) (1.7,3.6) (2.5,2.5)}--(1.2,5.3) -- (1,5.3) -- (1,1) -- (2.5,2.5)--cycle;
	\fill[domain=2:4.8,color=black!10!white] plot[smooth, tension=.7] coordinates {(2.5,2.5) (3.6,1.7) (5.3,1.2)}--(5.3,1.2) -- (5.3,1) -- (1,1) -- (2.5,2.5)--cycle;
	\fill[domain=2:4.8,color=black!20!white] plot[smooth, tension=.7] coordinates {(1.2,5.3) (1.7,3.6) (2.5,2.5)}--(2.5,2.5) -- (5.3,5.3) -- (1.2,5.3)--cycle;
	\fill[domain=2:4.8,color=black!20!white] plot[smooth, tension=.7] coordinates {(2.5,2.5) (3.6,1.7) (5.3,1.2)}--(5.3,1.2) -- (5.3,5.3) -- (2.5,2.5)--cycle;

	\draw[->] (-0.2,0) -- (5.8,0) node[below] {$p$};
	\draw[->] (0,-0.2) -- (0,5.4) node[left] {$q$};
	\draw[dashed, color=black]  (0, 0)--(5.3,5.3);
	\node[right] at (3.3,3.3) {{$\leftarrow$ $p=q$}};
	\node[left] at (0,-0.2) {{$0$}};
	\draw[fill] (1,0) circle[radius=1pt];
	\draw[fill] (2.5,2.5) circle[radius=1pt];
	\node[below] at (1,0) {{$1$}};
	\draw[fill] (0,1) circle[radius=1pt];
	\node[left] at (0,1) {{$1$}};
	\draw[fill] (2.5,0) circle[radius=1pt];
	\node[below] at (2.5,0) {{$p_{\mathrm{Fuj}}(n)$}};
	\draw[fill] (0,2.5) circle[radius=1pt];
	\node[left] at (0,2.5) {{$p_{\mathrm{Fuj}}(n)$}};
	\node[left, color=red] at (4.3,4.5) {{$\leftarrow$ $\alpha_{\max}(p,q)=\frac{n}{2}$}};
	\node[left, color=red] at (9,1.4) {{$\longleftarrow$ Blow-up in finite time}};
	\draw[dashed, color=black]  (0, 1)--(5.6, 1);
	\draw[dashed, color=black]  (1, 0)--(1, 5.4);
	\draw[dashed, color=black] (2.5,0)--(2.5,2.5);
	\draw[dashed, color=black] (0,2.5)--(2.5,2.5);
	\draw[color=red] plot[smooth, tension=.7] coordinates {(1.2,5.3) (1.7,3.6) (2.5,2.5)};
	\draw[color=red] plot[smooth, tension=.7] coordinates {(2.5,2.5) (3.6,1.7) (5.3,1.2)};

	\fill[color=black!20!white] (6,5)--(6.4,5)--(6.4,4.6)--(6,4.6)--cycle;
	\node[right] at (6.5,4.8) {{Global solution existence}};
	\fill[color=black!10!white] (6,4.4)--(6.4,4.4)--(6.4,4)--(6,4)--cycle;
	\node[right] at (6.5,4.2) {{Blow-up in finite time}};
	\end{tikzpicture}
	\caption{The critical curve for the coupled system \eqref{Eq_Coupled_Damped_Waves} in the $p-q$ plane}
	\label{imggg}
\end{figure}

\noindent To be specific, the authors in \cite{Sun-Wang-2007} obtained the critical condition \eqref{Critical_Condition} for $n=1,3$. More precisely, if $\alpha_{\max}(p,q)<n/2$, then there exists  a unique global (in time) small data Sobolev solution, whereas if $\alpha_{\max}(p,q)\geqslant n/2$, then every non-trivial local (in time) weak solution, in general, blows up in finite time. Afterwards, the same desired results for the dimensional cases $n=1,2,3$ were generalized in \cite{Nara-2009,Takeda-2009}, especially, some of decay estimates for solutions in time were improved for $n=3$. Finally, the recent papers \cite{Nish-Waka-2014,Nish-Waka-2015} demonstrated the critical condition \eqref{Critical_Condition} for any spatial dimensions $n\geqslant 1$, where the almost sharp estimates for the lifespan in the subcritical case $\alpha_{\max}(p,q)>n/2$ were found out, namely,
\begin{align*}
	\varepsilon^{-\frac{1}{\alpha_{\max}(p,q)-\frac{n}{2}}+\epsilon_0}\lesssim T_{\varepsilon}\lesssim \varepsilon^{-\frac{1}{\alpha_{\max}(p,q)-\frac{n}{2}}}
\end{align*} 
for any small number $\epsilon_0>0$. Hereafter the unexpressed multiplicative constants may depend on $n,p,q,u_0,u_1,v_0,v_1$ but are independent of $\varepsilon$. The authors also claimed some lower bound estimates for the lifespan in the critical case $\alpha_{\max}(p,q)=n/2$, nevertheless, it seems to be far from the (almost) sharp estimates. Again, we stress out that not only for the weakly coupled system of semilinear damped wave equations but also for the weakly coupled system of semilinear reaction-diffusion equations, both sharp upper bound and lower bound estimates of the lifespan in the critical case, i.e. under the critical condition \eqref{Critical_Condition}, are still completely open as far as the authors know. For this reason, our purpose of this paper is to give a positive answer, i.e. the sharp estimates \eqref{Sharp_Lifespan}, for these open problems.

To explore upper bound estimates for the lifespan of solutions to \eqref{Eq_Coupled_Damped_Waves} in the critical case \eqref{Critical_Condition}, we will employ the so-called \emph{test function method}, motivated by the recent studies \cite{Iked-Soba-2019,Dao-Reis-2020}. By constructing two test functions with suitable scaling (different from the scaling of the single semilinear damped wave equation in \cite{Iked-Soba-2019}), we will derive a system of two nonlinear differential inequalities with their initial values. Then, after dealing with these nonlinear inequalities, we will arrive at upper bound estimates for the lifespan with the help of some parameter-dependent auxiliary functionals.

To derive lower bound estimates for the lifespan of solutions to \eqref{Eq_Coupled_Damped_Waves}, we will introduce suitable Sobolev spaces with their corresponding norms carrying suitable polynomial-logarithmic type time-dependent weighted functions. Under this frame, employing the classical Gagliardo-Nirenberg inequality, we will get lower bound estimates for the lifespan in $n=1,2$, which are the same as the upper ones. From this observation, these obtained results are to conclude the sharpness of lifespan estimates immediately.


Since our approaches also can be applicable to the weakly coupled system of semilinear reaction-diffusion equations, we will propose some remarks in Section \ref{Section_Concluding_Remarks} on the sharp lifespan estimates for solutions in the critical case. Actually, it brings a possible answer for the question proposed in \cite{Fuji-Iked-Waka-2020}.

\medskip

\noindent\textbf{Notation:} Throughout this paper, we write $f\lesssim g$ when there exists a positive constant $C$ such that $f\leqslant Cg$. Moreover, $B_r$ stands for the ball around the origin with radius $r$ in $\mathbb{R}^n$.

\section{Main results}
To begin with this section, let us introduce some definitions of solutions to the semilinear Cauchy problem \eqref{Eq_Coupled_Damped_Waves}, which are to provide well-defined notations for the lifespan of corresponding solutions.
\begin{defn}\label{Defn_mild}
	The pair of functions $(u,v)$ is called a mild solution to the Cauchy problem \eqref{Eq_Coupled_Damped_Waves} on $[0,T)$ with $T>0$, if
	\begin{align}\label{Supp_01}
		(u,v)\in\ml{C}\big([0,T),H^1(\mb{R}^n)\big)\times \ml{C}\big([0,T),H^1(\mb{R}^n)\big)
	\end{align}
carrying its initial data satisfies the following integral systems:
\begin{align}\label{Represen_00}
	\begin{cases}
	\displaystyle{u(t,x)=\varepsilon(\partial_t+1)\ml{H}(t,\nabla)u_0(x)+\varepsilon\ml{H}(t,\nabla)u_1(x)+\int_0^t\ml{H}(t-\tau,\nabla)|v(\tau,x)|^p\mathrm{d}\tau,}\\
	\displaystyle{v(t,x)=\varepsilon(\partial_t+1)\ml{H}(t,\nabla)v_0(x)+\varepsilon\ml{H}(t,\nabla)v_1(x)+\int_0^t\ml{H}(t-\tau,\nabla)|u(\tau,x)|^q\mathrm{d}\tau,}
	\end{cases}
\end{align}
for any $t\in[0,T)$ with the operator
\begin{align*}
	\ml{H}(t,\nabla):=\mathrm{e}^{-\frac{t}{2}}\frac{\sin\left(t\sqrt{|\nabla|^2-1/4}\,\right)}{\sqrt{|\nabla|^2-1/4}}.
\end{align*}
It can be defined by the use of Fourier transforms
\begin{align*}
	\ml{H}(t,\nabla)=\ml{F}_{\xi\to x}^{-1}\left(\mathrm{e}^{-\frac{t}{2}}\frac{\sin\left(t\sqrt{|\xi|^2-1/4}\,\right)}{\sqrt{|\xi|^2-1/4}}\right).
\end{align*}
\end{defn}
The lifespan of a mild solution in the sense of Definition \ref{Defn_mild} is denoted by $T_{\varepsilon,\mathrm{m}}$, whose definition is similar to \eqref{Lifespan_Defn}.

\begin{defn}\label{Defn_Weak}
	The pair of functions $(u,v)$ is called a weak solution to the Cauchy problem \eqref{Eq_Coupled_Damped_Waves} on $[0,T)$ with $T>0$, if 
	\begin{align}\label{Supp_02}
		(u,v)\in L_{\mathrm{loc}}^q\big([0,T)\times\mb{R}^n\big)\times L_{\mathrm{loc}}^p\big([0,T)\times\mb{R}^n\big)
	\end{align}
satisfies the following integral equalities:
\begin{align}\label{Integral_01}
&\int_0^T\int_{\mb{R}^n}\left(\partial_t^2\Psi_1(t,x)-\Delta\Psi_1(t,x)-\partial_t\Psi_1(t,x)\right)u(t,x)\mathrm{d}x\mathrm{d}t\notag\\
&\qquad=\int_0^T\int_{\mb{R}^n}\Psi_1(t,x)|v(t,x)|^p\mathrm{d}x\mathrm{d}t+ \varepsilon\int_{\mb{R}^n}\big(\Psi_1(0,x)(u_0(x)+u_1(x))-\partial_t\Psi_1(0,x)u_0(x)\big)\mathrm{d}x,
\end{align}
as well as
\begin{align}\label{Integral_02}
	&\int_0^T\int_{\mb{R}^n}\left(\partial_t^2\Psi_2(t,x)-\Delta\Psi_2(t,x)-\partial_t\Psi_2(t,x)\right)v(t,x)\mathrm{d}x\mathrm{d}t\notag\\
	&\qquad=\int_0^T\int_{\mb{R}^n}\Psi_2(t,x)|u(t,x)|^q\mathrm{d}x\mathrm{d}t+ \varepsilon\int_{\mb{R}^n}\big(\Psi_2(0,x)(v_0(x)+v_1(x))-\partial_t\Psi_2(0,x)v_0(x)\big)\mathrm{d}x,
\end{align}
for any $\Psi_1,\Psi_2\in\ml{C}_0^{\infty}([0,T)\times\mb{R}^n)$.
\end{defn}
The lifespan of a weak solution in the sense of Definition \ref{Defn_Weak} is denoted by $T_{\varepsilon,\mathrm{w}}$, whose definition is similar to \eqref{Lifespan_Defn}.

\begin{remark} \label{Lifespan.Relation.Remark}
\fontshape{n}
\selectfont
We notice that if $(u,v)$ is a mild solution to \eqref{Eq_Coupled_Damped_Waves} in the sense of Definition \ref{Defn_mild}, then $(u,v)$ is also a weak solution to \eqref{Eq_Coupled_Damped_Waves} in the sense of Definition \ref{Defn_Weak}. This statement can be easily indicated by the standard density argument (see, for example, Proposition 3.1 in \cite{Iked-Waka-2013}). Therefore, the following relation is obviously true:
\begin{align} \label{Lifespan.Relation}
	T_{\varepsilon,\mathrm{m}}\leqslant T_{\varepsilon}\leqslant T_{\varepsilon,\mathrm{w}}.
\end{align}
For this reason, we want to underline in advance that the proof of sharp lifespan results comes from estimating $T_{\varepsilon,\mathrm{w}}$ from the above and $T_{\varepsilon,\mathrm{m}}$ from the below.
\end{remark}

Concerning upper bound estimates for the lifespan $T_{\varepsilon,\mathrm{w}}$ to \eqref{Eq_Coupled_Damped_Waves} in the critical case for all $n\geqslant 1$, we state the following result.
\begin{theorem}\label{Thm_Upper_Bound}
	Let us assume that initial data $u_j,v_j \in \mathcal{C}^\ity_0(\R^n)$ with $j=0,1$ satisfy
	\begin{align}\label{Assumption_Initial_data}
		\int_{\R^n}(u_0(x)+ u_1(x))\mathrm{d}x>0 \ \  \text{and}\ \  \int_{\R^n}(v_0(x)+ v_1(x))\mathrm{d}x>0.
	\end{align}
	If $p,q>1$ fulfill the critical condition \eqref{Critical_Condition} and $1<p,q\leqslant n/(n-2)$ if $n\geqslant 3$, then there exists a positive constant $\varepsilon_0$ such that for any $\e \in (0,\e_0]$ the lifespan $T_{\varepsilon,\mathrm{w}}$ of weak solutions to the Cauchy problem \eqref{Eq_Coupled_Damped_Waves} possesses the following upper bounds:
	\begin{align}\label{Lifespan}
		T_{\varepsilon,\mathrm{w}} \leqslant\begin{cases}
			\mathrm{exp}\left(C\e^{-(p-1)}\right) &\text{if}\ \ p=q, \\
			\mathrm{exp}\left(C\e^{-\max\left\{\frac{p(pq-1)}{p+1},\frac{q(pq-1)}{q+1}\right\}}\right) &\text{if}\ \ p\neq q,
		\end{cases}
	\end{align}
	where $C$ is a positive constant independent of $\e$.
\end{theorem}
\begin{remark} \label{remark2.1}
\fontshape{n}
\selectfont
Regarding the special case $p=q$ in our model \eqref{Eq_Coupled_Damped_Waves}, the critical condition \eqref{Critical_Condition} can be reduced to $p=q=p_{\mathrm{Fuj}}(n)$. Under this situation, our obtained result in \eqref{Lifespan} exactly coincides with the sharp upper bound estimate for the lifespan of solutions to the single semilinear damped wave equation \eqref{Eq_Single_Semilinear_Damped_Wave} in the critical case $p=p_{\mathrm{Fuj}}(n)$. Involving the latter issue, one may see \cite{Li-Zhou-1995,Lai-Zhou-2019,Iked-Soba-2019} for more details.
\end{remark}
\begin{remark}
\fontshape{n}
\selectfont
The condition $1<p,q\leqslant n/(n-2)$ if $n\geqslant 3$, appearing in Theorem \ref{Thm_Upper_Bound}, is to guarantee the local (in time) existence of solutions (see Proposition 2.1 in \cite{Nish-Waka-2015}).
\end{remark}
\begin{remark}
\fontshape{n}
\selectfont
Concerning the non-symmetric case $p\neq q$ in Theorem \ref{Thm_Upper_Bound}, from the critical condition \eqref{Critical_Condition} we get $\max\{p,q\}=n(pq-1)/2-1$, so that we may rewrite
	\begin{align*}
		\mathrm{exp}\left(C\e^{-\max\left\{\frac{p(pq-1)}{p+1},\frac{q(pq-1)}{q+1}\right\}}\right)= \mathrm{exp}\left(C\e^{-(pq-p_{\mathrm{Fuj}}(n))}\right).
	\end{align*}
In other words, it provides a way to see the lifespan estimates as in \eqref{Sharp_Lifespan}. This new discovery is one of the cores of this paper.
\end{remark}

To guarantee the sharpness of the derived lifespan estimates \eqref{Lifespan}, we have to estimate the lifespan $T_{\varepsilon,\mathrm{m}}$ from the below. Thus, we turn to lower bound estimates in the subsequent theorem.
\begin{theorem}\label{Thm_Lower_Bound}
	Let us assume that initial date belong to the following classical energy space with additional $L^1$ regularity:
	\begin{align*}
	\left((u_0,u_1),(v_0,v_1)\right) \in \mathcal{D}:= \left(\left(H^1(\R^n)\cap L^1(\R^n)\right) \times \left(L^2(\R^n)\cap L^1(\R^n)\right)\right)^2
	\end{align*}
	for $n=1,2$ with the corresponding norm
	\begin{align*}
		J[u_0,u_1,v_0,v_1]:= \left\|\left((u_0,u_1),(v_0,v_1)\right)\right\|_{\mathcal{D}} &= \|u_0\|_{H^1(\R^n)}+ \|u_0\|_{L^1(\R^n)}+ \|u_1\|_{L^2(\R^n)}+ \|u_1\|_{L^1(\R^n)} \\
		&\quad + \|v_0\|_{H^1(\R^n)}+ \|v_0\|_{L^1(\R^n)}+ \|v_1\|_{L^2(\R^n)}+ \|v_1\|_{L^1(\R^n)}.
	\end{align*}
	Moreover, we suppose that $p,q$ fulfill the critical condition \eqref{Critical_Condition} with $p,q\geqslant 2$ if $n=1,2$. Then, there exists a positive constant $\varepsilon_0$ such that for any $\e \in (0,\e_0]$ the lifespan $T_{\varepsilon,\mathrm{m}}$ of mild solutions to the Cauchy problem \eqref{Eq_Coupled_Damped_Waves} enjoys the following lower bounds:
	\begin{align}\label{Lifespan_Lower}
		T_{\varepsilon,\mathrm{m}} \geqslant\begin{cases}
			\mathrm{exp}\left(c\e^{-(p-1)}\right) &\text{if}\ \ p=q, \\
			\mathrm{exp}\left(c\e^{-\max\left\{\frac{p(pq-1)}{p+1},\frac{q(pq-1)}{q+1}\right\}}\right) &\text{if}\ \ p\neq q,
		\end{cases}
	\end{align}
	where $c$ is a positive constant depending on $n$ and $J[u_0,u_1,v_0,v_1]$ only.
\end{theorem}
\begin{remark}
\fontshape{n}
\selectfont
Similarly to Remark \ref{remark2.1}, our achieved result in \eqref{Lifespan_Lower} in the symmetric case $p=q=p_{\mathrm{Fuj}}(n)$ also exactly coincides with the sharp lower bound estimate for the lifespan of solutions to the single semilinear damped wave equation \eqref{Eq_Single_Semilinear_Damped_Wave} in the critical case $p=p_{\mathrm{Fuj}}(n)$.
\end{remark}
\begin{remark}
\fontshape{n}
\selectfont
Summarizing the derived results in Theorems \ref{Thm_Upper_Bound} and \ref{Thm_Lower_Bound} combined with the relation \eqref{Lifespan.Relation}, we claim that the sharp lifespan estimates $T_{\varepsilon}$ for solutions to the Cauchy problem \eqref{Eq_Coupled_Damped_Waves} in the critical case \eqref{Critical_Condition} are given by 
\begin{align*}
	T_{\varepsilon}\sim\begin{cases}
		\exp\left( C\varepsilon^{-(p-1)} \right)&\mbox{if}\ \ p=q,\\
		\exp\left(C\e^{-(pq-p_{\mathrm{Fuj}}(n))}\right) &\mbox{if}\ \ p\neq q,
	\end{cases}
\end{align*}
in low spatial dimensions, with a positive constant $C$ independent of $\varepsilon$. It seems also interesting to generalize these lower bound estimates for higher spatial dimensions by introducing weighted Sobolev spaces as well as employing some weighted decay estimates with respect to spatial variables. Namely, we conjecture that the sharp lifespan estimates above still hold for any $n\geqslant 1$. However, this purpose is beyond the scope of our paper.
\end{remark}

\section{Proof of Theorem \ref{Thm_Upper_Bound}}\label{Sec_Proof_Upper}

\subsection{Setting and test functions}\label{SubSec_Test_Fun}
Let us define the size of supports for initial data by
\begin{align*}
	r_0&:=\max\left\{ |x|: x\in\mathrm{supp}\, u_0\cup \mathrm{supp}\, u_1 \right\},\\
	r_1&:=\max\left\{ |x|: x\in\mathrm{supp}\, v_0\cup \mathrm{supp}\, v_1 \right\}.
\end{align*}
\begin{remark}
\fontshape{n}
\selectfont
We do not use finite propagation speed of solutions to damped wave equations. One recognizes that the support conditions of initial data will give remarkable contributions to catching the upper bound estimates for the lifespan of solutions. For this reason, it would provide an effective way to generalize our approach to some models without hyperbolic structure, for example, the weakly coupled system for reaction-diffusion equations \eqref{Eq_Coupled_Heats}.
\end{remark}
\noindent Without loss of generality, we  assume
\begin{align*}
	R_0:=\sqrt[4]{2\max\left\{r_0^4,r_1^4\right\}}<\sqrt{T_{\varepsilon,\mathrm{w}}}.
\end{align*}
We now introduce a test function $\eta=\eta(s)$ such that
\begin{align*}
\eta\in\ml{C}_0^{\infty}([0,\infty))\ \ \mbox{and}\ \ \eta(s):=\begin{cases}
	1&\mbox{if}\ \ s\in[0,1/2],\\
	\mbox{decreasing}&\mbox{if}\ \ s\in(1/2,1),\\
	0&\mbox{if}\ \ s\in[1,\infty).
\end{cases}
\end{align*}
Moreover, another test function $\eta^*=\eta^*(s)$ is also introduced by
\begin{align*}
 \eta^*(s):=\begin{cases}
		0&\mbox{if}\ \ s\in[0,1/2),\\
	\eta(s)&\mbox{if}\ \ s\in[1/2,\infty).
\end{cases}
\end{align*}
Next, for a large parameter $R\in(0,\infty)$, we take $\psi_R=\psi_R(t,x)$ and $\psi_R^*=\psi_R^*(t,x)$, which are defined, respectively, by
\begin{align*}
\psi_R(t,x):=\left(\eta\left(\frac{t^2+|x|^4}{R^4}\right)\right)^{\mu+2}\ \ \mbox{and}\ \ \psi_R^*(t,x):=\left(\eta^*\left(\frac{t^2+|x|^4}{R^4}\right)\right)^{\mu+2}
\end{align*}
with a positive constant $\mu$ fulfilling
\begin{align*}
	\mu\geqslant \max \left\{ \frac{2}{p-1},\frac{2}{q-1} \right\}.
\end{align*}
\subsection{Upper bound estimates for the lifespan}\label{SubSec_Upper_Bound}
By multiplying the test function $\psi_R$ on the both sides of the first and second equations in the Cauchy problem \eqref{Eq_Coupled_Damped_Waves} as well as integrating the resultants over $\mb{R}^n$, we obtain
\begin{align}\label{Eq_01}
\int_{\mb{R}^n}|v(t,x)|^p\psi_R(t,x)\mathrm{d}x&=\frac{\mathrm{d}^2}{\mathrm{d}t^2}\int_{\mb{R}^n}u(t,x)\psi_R(t,x)\mathrm{d}x+\frac{\mathrm{d}}{\mathrm{d}t}\int_{\mb{R}^n}\left(u(t,x)\psi_R(t,x)-2u(t,x)\partial_t\psi_R(t,x)\right)\mathrm{d}x\notag\\
&\quad+\int_{\mb{R}^n}u(t,x)\left(\partial_t^2\psi_R(t,x)-\Delta\psi_R(t,x)-\partial_t\psi_R(t,x)\right)\mathrm{d}x,
\end{align}
and similarly,
\begin{align*}
	\int_{\mb{R}^n}|u(t,x)|^q\psi_R(t,x)\mathrm{d}x&=\frac{\mathrm{d}^2}{\mathrm{d}t^2}\int_{\mb{R}^n}v(t,x)\psi_R(t,x)\mathrm{d}x+\frac{\mathrm{d}}{\mathrm{d}t}\int_{\mb{R}^n}\left(v(t,x)\psi_R(t,x)-2v(t,x)\partial_t\psi_R(t,x)\right)\mathrm{d}x\notag\\
	&\quad+\int_{\mb{R}^n}v(t,x)\left(\partial_t^2\psi_R(t,x)-\Delta\psi_R(t,x)-\partial_t\psi_R(t,x)\right)\mathrm{d}x,
\end{align*}
where we used integration by parts with respect to spatial variables and $\psi_R(t,x)\equiv0$ as $|x|\to\infty$ from the support condition. From straightforward computations, we observe
\begin{align*}
	\partial_t\psi_R(t,x)&=\frac{2(\mu+2)t}{R^4}\left(\eta\left(\frac{t^2+|x|^4}{R^4}\right)\right)^{\mu+1}\eta'\left(\frac{t^2+|x|^4}{R^4}\right),\\
	\partial_t^2\psi_R(t,x)&=\frac{2(\mu+2)}{R^4}\left(\eta\left(\frac{t^2+|x|^4}{R^4}\right)\right)^{\mu+1}\eta'\left(\frac{t^2+|x|^4}{R^4}\right)\\
	&\quad+\frac{4(\mu+1)(\mu+2)t^2}{R^8}\left(\eta\left(\frac{t^2+|x|^4}{R^4}\right)\right)^{\mu}\left(\eta'\left(\frac{t^2+|x|^4}{R^4}\right)\right)^2\\
	&\quad+\frac{4(\mu+2)t^2}{R^8}\left(\eta\left(\frac{t^2+|x|^4}{R^4}\right)\right)^{\mu+1}\eta''\left(\frac{t^2+|x|^4}{R^4}\right),
\end{align*}
and
\begin{align*}
	\partial_{x_k}^2\psi_R(t,x)&=\frac{4(\mu+2)(|x|^2+2x_k^2)}{R^4}\left(\eta\left(\frac{t^2+|x|^4}{R^4}\right)\right)^{\mu+1}\eta'\left(\frac{t^2+|x|^4}{R^4}\right)\\
	&\quad+\frac{16(\mu+1)(\mu+2)|x|^4x_k^2}{R^8}\left(\eta\left(\frac{t^2+|x|^4}{R^4}\right)\right)^{\mu}\left(\eta'\left(\frac{t^2+|x|^4}{R^4}\right)\right)^2\\
	&\quad+\frac{16(\mu+2)|x|^4x_k^2}{R^8}\left(\eta\left(\frac{t^2+|x|^4}{R^4}\right)\right)^{\mu+1}\eta''\left(\frac{t^2+|x|^4}{R^4}\right).
\end{align*}
Due to the fact that
\begin{align*}
	\eta'\left(\frac{t^2+|x|^4}{R^4}\right)\not\equiv0,\ \ \eta''\left(\frac{t^2+|x|^4}{R^4}\right)\not\equiv0 \ \ \mbox{for}\ \ \frac{R^4}{2}<t^2+|x|^4<R^4,
\end{align*}
as well as $\eta\in\ml{C}_0^{\infty}([0,\infty))$, we are able to state
\begin{align*}
	\left|\partial_t^2\psi_R(t,x)-\Delta\psi_R(t,x)-\partial_t\psi_R(t,x)\right|&\lesssim\frac{1}{R^2}(\psi_R^*(t,x))^{\frac{\mu+1}{\mu+2}}+\frac{R^2+1}{R^4}(\psi_R^*(t,x))^{\frac{\mu}{\mu+2}}\\
	&\lesssim \frac{1}{R^2}(\psi_R^*(t,x))^{\frac{\mu}{\mu+2}}.
\end{align*}
In the above estimate, we have utilized $0<\psi_R^*(t,x)< 1$ and $R\gg1$. Taking account of $R\in[R_0,\sqrt{T_{\varepsilon,\mathrm{w}}})$ and integrating \eqref{Eq_01} over $(0,T_{\varepsilon,\mathrm{w}})$, we may deduce
\begin{align*}
\int_0^{T_{\varepsilon,\mathrm{w}}}\int_{\mb{R}^n}|v(t,x)|^p\psi_R(t,x)\mathrm{d}x\mathrm{d}t&\leqslant \left(\int_{\mb{R}^n}\left(u_t(t,x)\psi_R(t,x)+u(t,x)\partial_t\psi_R(t,x)\right)\mathrm{d}x\right)\Big|_{t=0}^{t=T_{\varepsilon,\mathrm{w}}}\\
&\quad+\left(\int_{\mb{R}^n}\left(u(t,x)\psi_R(t,x)-2u(t,x)\partial_t\psi_R(t,x)\right)\mathrm{d}x\right)\Big|_{t=0}^{t=T_{\varepsilon,\mathrm{w}}}\\
&\quad+\frac{1}{R^2}\int_0^{T_{\varepsilon,\mathrm{w}}}\int_{\mb{R}^n}|u(t,x)|(\psi_R^*(t,x))^{\frac{\mu}{\mu+2}}\mathrm{d}x\mathrm{d}t.
\end{align*}
The consideration $R\in[R_0,\sqrt{T_{\varepsilon,\mathrm{w}}})$ leads to
\begin{align*}
	2\max\{ r_0^4,r_1^4 \}=R_0^4\leqslant R^4\leqslant (T_{\varepsilon,\mathrm{w}})^2,
\end{align*}
so that
\begin{align*}
\psi_R(T_{\varepsilon,\mathrm{w}},x)=\partial_t\psi_R(T_{\varepsilon,\mathrm{w}},x)=0\ \ &\mbox{for any}\ \ x\in\mb{R}^n,\\
\psi_R(0,x)=1\ \ &\mbox{for any}\ \ x\in B_{r_0}.
\end{align*}
In other words, we have
\begin{align}\label{Eq_03}
&\varepsilon\underbrace{\int_{\mb{R}^n}(u_0(x)+u_1(x))\mathrm{d}x}_{=:I_0[u_0,u_1]}+\int_0^{T_{\varepsilon,\mathrm{w}}}\int_{\mb{R}^n}|v(t,x)|^p\psi_R(t,x)\mathrm{d}x\mathrm{d}t\notag\\
&\qquad \leqslant\frac{1}{R^2}\int_0^{T_{\varepsilon,\mathrm{w}}}\int_{\mb{R}^n}|u(t,x)|(\psi_R^*(t,x))^{\frac{\mu}{\mu+2}}\mathrm{d}x\mathrm{d}t\notag\\
&\qquad \leqslant\frac{1}{R^2}\left(\int_{\mathrm{supp}\,\psi_R^*}\mathrm{d}(x,t)\right)^{\frac{1}{q'}}\left(\int_0^{T_{\varepsilon,\mathrm{w}}}\int_{\mb{R}^n}|u(t,x)|^q(\psi_R^*(t,x))^{\frac{q\mu}{\mu+2}}\mathrm{d}x\mathrm{d}t\right)^{\frac{1}{q}}\notag\\
&\qquad \lesssim R^{n-\frac{n+2}{q}}\left(\int_0^{T_{\varepsilon,\mathrm{w}}}\int_{\mb{R}^n}|u(t,x)|^q(\psi_R^*(t,x))^{\frac{q\mu}{\mu+2}}\mathrm{d}x\mathrm{d}t\right)^{\frac{1}{q}},
\end{align}
where the support conditions for initial data since $\mathrm{supp}\,u_0\cup\mathrm{supp}\,u_1\subset B_{r_0}$ were used as well as the following attention should be recognized:
$$ \mathrm{supp}\,\psi^*_R \subset \left([0,R^2] \times B_R\right) \backslash \left\{(t,x) \,:\, t^2+|x|^4 \leqslant\frac{R^4}{2}\right\}. $$
Repeating the same procedure as the above, it holds
\begin{align}\label{Eq_04}
	&\varepsilon\underbrace{\int_{\mb{R}^n}(v_0(x)+v_1(x))\mathrm{d}x}_{=:I_0[v_0,v_1]}+\int_0^{T_{\varepsilon,\mathrm{w}}}\int_{\mb{R}^n}|u(t,x)|^q\psi_R(t,x)\mathrm{d}x\mathrm{d}t\notag\\	
	&\qquad \lesssim R^{n-\frac{n+2}{p}}\left(\int_0^{T_{\varepsilon,\mathrm{w}}}\int_{\mb{R}^n}|v(t,x)|^p(\psi_R^*(t,x))^{\frac{p\mu}{\mu+2}}\mathrm{d}x\mathrm{d}t\right)^{\frac{1}{p}}.
\end{align}
Let us introduce two auxiliary functionals as follows:
\begin{align*}
Y_q(R):=\int_0^Ry_q(r)r^{-1}\mathrm{d}r\ \ \mbox{with}\ \ y_q(r):=\int_0^{T_{\varepsilon,\mathrm{w}}}\int_{\mb{R}^n}|u(t,x)|^q(\psi_r^*(t,x))^{\frac{q\mu}{\mu+2}}\mathrm{d}x\mathrm{d}t,\\
Y_p(R):=\int_0^Ry_p(r)r^{-1}\mathrm{d}r\ \ \mbox{with}\ \ y_p(r):=\int_0^{T_{\varepsilon,\mathrm{w}}}\int_{\mb{R}^n}|v(t,x)|^p(\psi_r^*(t,x))^{\frac{p\mu}{\mu+2}}\mathrm{d}x\mathrm{d}t.
\end{align*}
The change of variable $s=(t^2+|x|^4)/r^4$  yields
\begin{align}
Y_q(R)&=\int_0^R \left(\int_0^{T_{\varepsilon,\mathrm{w}}}\int_{\mb{R}^n}|u(t,x)|^q(\psi_r^*(t,x))^{\frac{q\mu}{\mu+2}}\mathrm{d}x\mathrm{d}t\right)\,r^{-1}\mathrm{d}r \nonumber \\
&=\frac{1}{4}\int_0^{T_{\varepsilon,\mathrm{w}}}\int_{\mb{R}^n}|u(t,x)|^q\int_{(t^2+|x|^4)/R^4}^{\infty}(\eta^*(s))^{q\mu}s^{-1}\mathrm{d}s\mathrm{d}x\mathrm{d}t \nonumber \\
&\leqslant \frac{1}{4} \int_0^{T_{\varepsilon,\mathrm{w}}}\int_{\mb{R}^n}|u(t,x)|^q\left(\int_{1/2}^1 (\eta^*(s))^{\mu q}s^{-1}\mathrm{d}s\right)\mathrm{d}x\mathrm{d}t, \label{Eq_05}
\end{align}
where we considered the support condition for $\eta^*(s)$ in the third line of the chain estimates above. Thus, one may arrive at
\begin{align}
Y_q(R)&\leqslant\frac{1}{4} \int_0^{T_{\varepsilon,\mathrm{w}}}\int_{\mb{R}^n}|u(t,x)|^q\sup_{r \in (0,R)}\left(\eta^*\left(\frac{t^2+|x|^4}{r^4}\right)\right)^{\mu q} \left(\int_{1/2}^1 s^{-1}\mathrm{d}s\right)\mathrm{d}x\mathrm{d}t \nonumber \\
&\leqslant\frac{\log 2}{4} \int_0^{T_{\varepsilon,\mathrm{w}}}\int_{\mb{R}^n}|u(t,x)|^q \left(\eta^*\left(\frac{t^2+|x|^4}{R^4}\right)\right)^{\mu q}\mathrm{d}x\mathrm{d}t \nonumber \\
&\leqslant\frac{\log 2}{4}y_q(R). \label{Eq_06}
\end{align}
In addition, using the property $\eta^*(s)\equiv \eta(s)$ for any $s\in[1/2,1]$ in \eqref{Eq_05} we also verify the following estimate:
\begin{align}
Y_q(R) &\lesssim \int_0^{T_{\varepsilon,\mathrm{w}}}\int_{\mb{R}^n}|u(t,x)|^q \left(\eta\left(\frac{t^2+|x|^4}{R^4}\right)\right)^{\mu q}\mathrm{d}x\mathrm{d}t, \nonumber \\
&\lesssim \int_0^{T_{\varepsilon,\mathrm{w}}}\int_{\mb{R}^n}|u(t,x)|^q\left(\psi_R(t,x)\right)^{\frac{q\mu}{\mu+2}}\mathrm{d}x\mathrm{d}t. \label{Eq_07}
\end{align}
By the similar fashion, one gets
\begin{align}
Y_p(R)&\leqslant\frac{\log 2}{4}y_p(R), \label{Eq_08} \\
Y_p(R)&\lesssim \int_0^{T_{\varepsilon,\mathrm{w}}}\int_{\mb{R}^n}|v(t,x)|^p\left(\psi_R(t,x)\right)^{\frac{p\mu}{\mu+2}}\mathrm{d}x\mathrm{d}t. \label{Eq_09}
\end{align}
To derive adaptable functionals, we recall
\begin{align*}
	\mu\geqslant \max \left\{ \frac{2}{p-1},\frac{2}{q-1} \right\}=\frac{2}{\min\{p,q\}-1}, \ \ \text{i.e.}\ \  \frac{\min\{p,q\}\,\mu}{\mu +2}\geqslant 1
\end{align*}
to show from \eqref{Eq_07} and \eqref{Eq_09} that
\begin{align}
	Y_q(R)&\lesssim \int_0^{T_{\varepsilon,\mathrm{w}}}\int_{\mb{R}^n}|u(t,x)|^q\psi_R(t,x)\mathrm{d}x\mathrm{d}t,\label{Eq_10}\\
	Y_p(R)&\lesssim \int_0^{T_{\varepsilon,\mathrm{w}}}\int_{\mb{R}^n}|v(t,x)|^p\psi_R(t,x)\mathrm{d}x\mathrm{d}t.\label{Eq_11}
\end{align}
In what follows, we denote by $C_j$ with $j\in \N$ positive constants independent of $R$ and $\varepsilon$. According to the estimates \eqref{Eq_03}, \eqref{Eq_04}, \eqref{Eq_10}, \eqref{Eq_11} and the facts that
\begin{align*}
y_p(R)= RY'_p(R),\ \ y_q(R)= RY'_q(R),
\end{align*}
we conclude the following coupled system of nonlinear differential inequalities:
\begin{align}
	Y'_p(R)&\geqslant C_1\delta R^{n+1-np}\left(Y_q(R)+\varepsilon I[v_0,v_1]\right)^p,\label{Pre_Frame_1}\\
	Y'_q(R)&\geqslant C_2\delta R^{n+1-nq}\left(Y_p(R)+\varepsilon I[u_0,u_1]\right)^q,\label{Pre_Frame_2}
\end{align}
for any $R\in[R_0,\sqrt{T_{\varepsilon,\mathrm{w}}})$. Here, a constant $\delta \in (0,1]$ will be chosen later.  By considering our assumption on initial data \eqref{Assumption_Initial_data}, we have derived the differential frames
\begin{align*}
	Y'_p(R)&\geqslant C_1\delta R^{n+1-np}(Y_q(R))^p,\\
	Y'_q(R)&\geqslant C_2\delta R^{n+1-nq}(Y_p(R))^q,
\end{align*}
for any $R\in[R_0,\sqrt{T_{\varepsilon,\mathrm{w}}})$, with their initial values (from the integration of the inequalities \eqref{Pre_Frame_1} and \eqref{Pre_Frame_2} over $[R_0,R]$, respectively)
\begin{align}
Y_p(R)&\geqslant \varepsilon^p C_1\delta (I[v_0,v_1])^p\int_{R_0}^{R}\rho^{n+1-np}\mathrm{d}\rho, \nonumber \\
Y_q(R)&\geqslant \varepsilon^q C_2\delta (I[u_0,u_1])^q\int_{R_0}^{R}\rho^{n+1-nq}\mathrm{d}\rho, \label{Eq_*}
\end{align}
under the restriction \eqref{Critical_Condition}.

Without loss of generality, we will focus on the treatment of the case $\max\{p,q\}=q\neq p$ only due to the fact is that the case $\max\{p,q\}=p\neq q$ can be also treated in the same way. The rest case $p=q$ will be shown later. Namely, to get started, the condition \eqref{Critical_Condition} can be written by
\begin{equation}
\frac{q+1}{pq-1}=\frac{n}{2}\ \ \mbox{and}\ \ q>1+\frac{2}{n}. \label{Eq_**}
\end{equation}
Let us now set up the two auxiliary functions $ \phi_1(R):= R^{n+1-np}$ and $\phi_2(R):= R^{n+1-nq}$ to re-express the differential frames above in the following way:
\begin{align}
Y_p'(R) &\geqslant C_1\delta\phi_1(R) (Y_q(R))^{p}, \label{Eq_12} \\
Y_q'(R) &\geqslant C_2\delta\phi_2(R) (Y_p(R))^{q}. \label{Eq_13}
\end{align}
Multiplying \eqref{Eq_12} by $Y_q'(R)$ and then carrying out integration by parts over $[R_0,R]$ give
\begin{align*}
&Y_p(R)Y_q'(R)- Y_p(R_0)Y_q'(R_0)- \int_{R_0}^R Y_p(s)Y_q''(s)\mathrm{d}s \\
&\qquad \geqslant \frac{C_1\delta}{p+1} \phi_1(R) (Y_q(R))^{p+1}- \frac{C_1\delta}{p+1} \phi_1(R_0) (Y_q(R_0))^{p+1} - \frac{C_1\delta}{p+1} \int_{R_0}^R \phi'_1(s) (Y_q(s))^{p+1}\mathrm{d}s.
\end{align*}
By the aid of the relation 
\begin{align}\label{Yq''}
Y_q''(s)= \frac{y_q'(s)- Y_q'(s)}{s},	
\end{align}
 one gains
\begin{align}
&Y_p(R)Y_q'(R)+ \int_{R_0}^R \frac{Y_p(s)Y_q'(s)}{s}\mathrm{d}s- \int_{R_0}^R \frac{Y_p(s)y_q'(s)}{s}\mathrm{d}s \nonumber \\
&\qquad \geqslant \frac{C_1\delta}{p+1} \phi_1(R) (Y_q(R))^{p+1}+ \left(Y_p(R_0)Y_q'(R_0)- \frac{C_1\delta}{p+1} \phi_1(R_0) (Y_q(R_0))^{p+1}\right) \nonumber \\
&\qquad \quad - \frac{C_1\delta}{p+1} \int_{R_0}^R \phi'_1(s) (Y_q(s))^{p+1}\mathrm{d}s. \label{Eq_14}
\end{align}
Let us devote our consideration to the estimate for the right-hand side (RHS) of \eqref{Eq_14}. In the first stage, we need to opt a constant $\delta=\delta(C_1,R_0,p)$ fulfilling
\begin{align}\label{Choice_delta}
	0< \delta\leqslant\min\left\{\frac{(p+1)Y_p(R_0)Y_q'(R_0)}{C_1\phi_1(R_0)(Y_q(R_0))^{p+1}},1\right\}
\end{align}
so that we may conclude that
\begin{align*}
Y_p(R_0)Y_q'(R_0)- \frac{C_1\delta}{p+1} \phi_1(R_0) (Y_q(R_0))^{p+1} \geqslant 0.	
\end{align*}
Noticing that $Y'_q(R_0)=y_q(R_0)/R_0>0$, the range of $\delta$ in \eqref{Choice_delta} is not empty. Subsequently, it implies
\begin{align}
\text{RHS of }\eqref{Eq_14} &\geqslant \frac{C_1\delta}{p+1} \phi_1(R) (Y_q(R))^{p+1}- \frac{C_1\delta}{p+1} \int_{R_0}^R \phi'_1(s) (Y_q(s))^{p+1}\mathrm{d}s \nonumber \\ 
&= \frac{C_1\delta}{p+1} \phi_1(R) (Y_q(R))^{p+1}- \frac{C_1\delta}{p+1}(n+1-np) \int_{R_0}^R s^{n-np} (Y_q(s))^{p+1}\mathrm{d}s. \label{Eq_15}
\end{align}
Concerning the sign for the last term in the previous inequality, our next arguments are divided into two cases separately as follows.
\begin{itemize}
\item \textbf{Case 1}: When $1+ \frac{1}{n}\leqslant p<1+\frac{2}{n}$, we have $n+1-np\leqslant0$. Then, by \eqref{Eq_15} it is obvious to catch the estimate
\begin{equation*}
\text{RHS of }\eqref{Eq_14}\geqslant \frac{C_1\delta}{p+1} \phi_1(R) (Y_q(R))^{p+1}.
\end{equation*}
\item \textbf{Case 2}: When $1< p< 1+ \frac{1}{n}$, we get $n+1-np> 0$. Then, setting $h_1=h_1(s)$ in the integrand of (\ref{Eq_15}) by
$$ h_1(s):= s^{n-np} (Y_q(s))^{p+1}, $$
we can calculate straightforwardly in this way
\begin{align*}
h'_1(s) &= (n-np) s^{n-np-1} (Y_q(s))^{p+1}+ (p+1)s^{n-np} (Y_q(s))^{p}Y'_q(s) \\
&= s^{n-np-1} (Y_q(s))^{p}\left((n-np)Y_q(s)+ (p+1)sY'_q(s)\right) \\
&= s^{n-np-1} (Y_q(s))^{p}\left((n-np)Y_q(s)+ (p+1)y_q(s)\right),
\end{align*}
where we noticed again that the relation $y_q(s)= sY'_q(s)$ holds. By using the derived estimate \eqref{Eq_06}, we can proceed as follows:
\begin{align*}
h'_1(s) &\geqslant \left((n-np)\frac{\log 2}{4}+ (p+1)\right)s^{n-np-1} (Y_q(s))^{p}y_q(s) \\ 
&> \left(p+1- \frac{\log 2}{4}\right)s^{n-np-1} (Y_q(s))^{p}y_q(s) \geqslant 0
\end{align*}
due to the strict inequality $n-np> -1$. As a result, $h_1= h_1(s)$ is a strictly increasing function so that we derive
\begin{equation}
\int_{R_0}^R s^{n-np} (Y_q(s))^{p+1}\mathrm{d}s \leqslant R^{n-np} (Y_q(R))^{p+1}(R-R_0) \leqslant\phi_1(R) (Y_q(R))^{p+1}. \label{Eq_17}
\end{equation}
Both the estimates \eqref{Eq_15} and \eqref{Eq_17} lead to
\begin{equation}
\text{RHS of }\eqref{Eq_14}\geqslant \frac{C_1\delta n(p- 1)}{p+1} \phi_1(R) (Y_q(R))^{p+1}. \label{Eq_18}
\end{equation}
\end{itemize}
Let us now come back to control the left-hand side (LHS) of \eqref{Eq_14}. At first, thanks to the non-decreasing property of $y_q= y_q(s)$, indeed,
\begin{align*}
	y_q'(s)=-\frac{4q\mu}{s^5}\int_0^{T_{\varepsilon,\mathrm{w}}}\int_{\mb{R}^n}|u(t,x)|^q\left(\eta^*\left(\frac{t^2+|x|^4}{s^4}\right)\right)^{q\mu-1}\left(\eta^*\left(\frac{t^2+|x|^4}{s^4}\right)\right)'(t^2+|x|^4)\mathrm{d}x\mathrm{d}t\geqslant0 
\end{align*}
with the help of non-increasing property of $\eta^*$, for any $s \in [R_0,R]$, one achieves
\begin{equation}
\text{LHS of }\eqref{Eq_14} \leqslant Y_p(R)Y_q'(R)+ \int_{R_0}^R \frac{Y_p(s)Y_q'(s)}{s}\mathrm{d}s. \label{Eq_19}
\end{equation}
Putting $h_2=h_2(s)$ in the integrand of the last term in \eqref{Eq_19} by
$$ h_2(s):= \frac{Y_p(s)Y_q'(s)}{s} $$
and noticing the equality \eqref{Yq''} we show that
\begin{align*}
h'_2(s) &= \frac{Y'_p(s)Y_q'(s)s+ Y_p(s)Y''_q(s)s- Y_p(s)Y_q'(s)}{s^2} \\
&= \frac{Y'_p(s)Y_q'(s)s+ Y_p(s) \left(y'_q(s)- Y'_q(s)\right)- Y_p(s)Y_q'(s)}{s^2}.
\end{align*}
Moreover, thanks to $y'_q(s)\geqslant 0$ and $y_p(s)=sY'_p(s)$, we deduce
\begin{align*}
h'_2(s) \geqslant \frac{Y_q'(s)\left(Y'_p(s)s- 2Y_p(s)\right)}{s^2} &= \frac{Y_q'(s)\big(y_p(s)- 2Y_p(s)\big)}{s^2} \\
&\geqslant \left(1- \frac{\log 2}{2}\right)\frac{y_p(s)Y_q'(s)}{s^2} \geqslant 0,
\end{align*}
where we have employed the obtained inequality \eqref{Eq_08}. This means that $h_2= h_2(s)$ is a non-decreasing function. In other words, it follows
\begin{equation}
\int_{R_0}^R \frac{Y_p(s)Y_q'(s)}{s}\mathrm{d}s \leqslant\frac{Y_p(R)Y_q'(R)}{R}(R-R_0)\leqslant Y_p(R)Y_q'(R). \label{Eq_20}
\end{equation}
For this reason, one may combine \eqref{Eq_19} and \eqref{Eq_20} to get
\begin{equation}
\text{LHS of }\eqref{Eq_14} \leqslant2Y_p(R)Y_q'(R). \label{Eq_21}
\end{equation}
Summarizing, the link of these derived estimates \eqref{Eq_14}, \eqref{Eq_18} and \eqref{Eq_21} is to indicate that
$$ Y_p(R)Y_q'(R) \geqslant C_0 \phi_1(R) (Y_q(R))^{p+1}, $$
which is equivalent to
\begin{equation}
Y_p(R) \geqslant \frac{C_0 \phi_1(R) (Y_q(R))^{p+1}}{Y_q'(R)} \label{Eq_22}
\end{equation}
for $R\geqslant R_0$. Hence, substituting \eqref{Eq_22} into \eqref{Eq_13} entails
$$ Y_q'(R) \geqslant \frac{C_3 (\phi_1(R))^{q}\phi_2(R)\,(Y_q(R))^{q(p+1)}}{(Y_q'(R))^{q}}, $$
which implies immediately
\begin{align*}
Y_q'(R) &\geqslant C_3^{\frac{1}{q+1}} (\phi_1(R))^{\frac{q}{q+1}}(\phi_2(R))^{\frac{1}{q+1}}(Y_q(R))^{\frac{q(p+1)}{q+1}} \\
& = C_3^{\frac{1}{q+1}}R^{1-\frac{n(pq-1)}{q+1}} (Y_q(R))^{\frac{q(p+1)}{q+1}}\\
&= C_3^{\frac{1}{q+1}}R^{-1} (Y_q(R))^{\frac{q(p+1)}{q+1}}
\end{align*}
in our case \eqref{Eq_**}. Clearly, the above estimate is to verify the following:
\begin{equation}
\frac{Y_q'(R)}{(Y_q(R))^{\frac{q(p+1)}{q+1}}} \geqslant C_3^{\frac{1}{q+1}}R^{-1}.\label{Eq_23}
\end{equation}
Then, considering $R\geqslant R_0^2$ we take integration of two sides of \eqref{Eq_23} over $[\sqrt{R},R]$ to obtain
\begin{align*}
-\frac{q+1}{pq-1} (Y_q(s))^{-\frac{pq-1}{q+1}}\Big|_{s=\sqrt{R}}^{s=R} &= \frac{n}{2}\left((Y_q(\sqrt{R}\,))^{-\frac{2}{n}}- (Y_q(R))^{-\frac{2}{n}}\right) \\ 
&\geqslant C_3^{\frac{1}{q+1}}\left(\log R- \log(\sqrt{R}\,)\right)= \frac{1}{2}C_3^{\frac{1}{q+1}}\log R.
\end{align*}
Therefore, it holds
\begin{align}\label{Eq_16}
	\log R\leqslant n C_3^{-\frac{1}{q+1}} \left(Y_q(\sqrt{R}\,)\right)^{-\frac{2}{n}}.
\end{align}
By recalling the inequality \eqref{Eq_*}, it is obvious to catch the estimate
$$ Y_q(\sqrt{R}\,) \geqslant \varepsilon^q C_4\delta (I[u_0,u_1])^q $$
since we are in the situation $n+1-nq<-1$ from \eqref{Eq_**} as well as $R\geqslant R_0$ to ensure the boundedness of the integral in \eqref{Eq_*}. Hence, one arrives at the next estimate by the combination of the last two inequalities
$$\log \sqrt{T_{\varepsilon,\mathrm{w}}}=\lim\limits_{R\uparrow \sqrt{T_{\varepsilon,\mathrm{w}}}} \log R\leqslant C_5\varepsilon^{-\frac{2q}{n}}= C_5\varepsilon^{-\frac{q(pq-1)}{q+1}}, $$
where we note that $\frac{2}{n}= \frac{pq-1}{q+1}$. This is to show the desired upper bound of lifespan estimate for mild solutions to the Cauchy problem \eqref{Eq_Coupled_Damped_Waves}.

We remark that in the special case $p=q=p_{\mathrm{Fuj}}(n)$, one finds from \eqref{Eq_*} the following estimate:
\begin{align*}
	Y_q(R)\geqslant C_6\varepsilon^{p_{\mathrm{Fuj}}(n)}\log R.
\end{align*}
Then, we combine the previous inequality with \eqref{Eq_16} to obtain
\begin{align*}
	\log \sqrt{T_{\varepsilon,\mathrm{w}}}=\lim\limits_{R\uparrow \sqrt{T_{\varepsilon,\mathrm{w}}}}\log R\leqslant C_7\varepsilon^{-\frac{2}{n}}=C_7\varepsilon^{-(p-1)}.
\end{align*}
Finally, taking the action of the exponential function gives the completeness of our proof.

\section{Proof of Theorem \ref{Thm_Lower_Bound}}\label{Sec_Proof_Lower}
\subsection{Philosophy of our approach}
With the same reason of the proof of Theorem \ref{Thm_Upper_Bound}, we are going to focus on the case $p< q$ only, and we will give some remarks for the special case $p=q$. For the sake of brevity, we put
\begin{align*}
	\gamma(p,q):= \frac{q-p}{pq-1}>0\ \ \mbox{and}\ \ \alpha(p,q):=\frac{n(p-1)}{2q}>0
\end{align*}
because of the hypothesis $p<q$ and $p>1$, respectively.

First of all, we introduce the evolution spaces $Y_1(T)$ and $Y_2(T)$ as follows:
\begin{align*}
	Y_j(T)=\ml{C}([0,T],H^1(\mb{R}^n))\ \ \mbox{for}\ \ j=1,2,
\end{align*}
carrying their corresponding norms
\begin{align*}
	\|u\|_{Y_1(T)}:=\sup\limits_{t\in[0,T]}\left((1+t)^{-\gamma(p,q)} (\log(\mathrm{e}+t))^{\alpha(p,q)}\ml{M}[u](t)\right)\ \ \mbox{and}\ \ \|v\|_{Y_2(T)}:=\sup\limits_{t\in[0,T]}\left(\ml{M}[v](t)\right),
\end{align*}
where we define
\begin{align*}
	\ml{M}[w](t):=(1+t)^{\frac{n}{4}}\|w(t,\cdot)\|_{L^2(\mb{R}^n)}+(1+t)^{\frac{n}{4}+\frac{1}{2}}\|\nabla w(t,\cdot)\|_{L^2(\mb{R}^n)},
\end{align*}
with $w=u$ or $w=v$. With the last definitions, we can introduce the solution space $X(T)$ of the weakly coupled system \eqref{Eq_Coupled_Damped_Waves} by
$$ X(T)= Y_1(T) \times Y_2(T), $$
endowed with the norm
$$ \|(u,v)\|_{X(T)}:=\|u\|_{Y_1(T)}+ \|v\|_{Y_2(T)}. $$
Because different power nonlinearities have different influence on lifespan estimates, we can allow the effect of the loss of decay in comparison with the corresponding homogeneous linear problem (see Proposition \ref{Estimates_Linear_Damped_Waves} later). Particularly, in our consideration $p<p_{\mathrm{Fuj}}(n)<q$, we take a loss $(1+t)^{-\gamma(p,q)}(\log(\mathrm{e+t}))^{\alpha(p,q)}$ in the norm of $\|u\|_{Y_1(T)}$. As we can see in Section \ref{Sec.4.2}, this polynomial-logarithmic loss of decay plays an essential role in our proof.

Let us now denote by $\mathcal{H}_0(t,x)$ and $\mathcal{H}_1(t,x)$, the fundamental solutions to the following linearized Cauchy problem:
\begin{align} \label{Eq_Linear_Damped_Waves}
	\begin{cases}
		w_{tt}-\Delta w+w_t= 0,&x\in\mb{R}^n,\ t\in(0,T),\\
		(w,w_t)(0,x)=( w_0, w_1)(x),&x\in\mb{R}^n,
	\end{cases}
\end{align}
for $w=u$ or $w=v$, so that we can represent the solution formula to the corresponding homogeneous linear Cauchy problem for \eqref{Eq_Coupled_Damped_Waves} by
$$
\begin{cases}
u^{\lin}(t,x):= \varepsilon \mathcal{H}_0(t,x)\ast_{(x)} u_{0}(x)+ \varepsilon \mathcal{H}_1(t,x)\ast_{(x)} u_{1}(x),\\
v^{\lin}(t,x):= \varepsilon \mathcal{H}_0(t,x)\ast_{(x)} v_{0}(x)+ \varepsilon \mathcal{H}_1(t,x)\ast_{(x)} v_{1}(x),
\end{cases}
$$
where $\ast_{(x)}$ stands for the convolution with respect to spatial variables $x$. Then, heavily motivated by Duhamel's principle, the solution to \eqref{Eq_Coupled_Damped_Waves} can be written in this form
$$
\begin{cases}
u(t,x)=u^{\lin}(t,x)+\displaystyle\int_{0}^{t}\mathcal{H}_1(t-\tau,x)\ast_{(x)} |v(\tau,x)|^{p}\mathrm{d}\tau=:u^{\lin}(t,x)+u^{\non}(t,x),
\\
v(t,x)=v^{\lin}(t,x)+\displaystyle\int_{0}^{t} \mathcal{H}_1(t-\tau,x)\ast_{(x)} |u(\tau,x)|^{q}\mathrm{d}\tau=:v^{\lin}(t,x)+v^{\non}(t,x),
\end{cases}
$$ where is an equivalent way to represent the mild solution \eqref{Represen_00}.

The main point of our approach to indicate the desired lower bound estimates for the lifespan relies on the proof of a pair of inequalities as follows:
\begin{align}
\|u\|_{Y_1(T)} &\leqslant\varepsilon c_0+ c^u_1 (\log(\mathrm{e}+t))^{\alpha(p,q) }\,\|v\|_{Y_2(T)}^{p}, \label{Ineq.01} \\
\|v\|_{Y_2(T)} &\leqslant\varepsilon c_0+ c^v_1 (\log(\mathrm{e}+t))^{1-\alpha(p,q) q}\,\|u\|_{Y_1(T)}^{q}, \label{Ineq.02}
\end{align}
for all $t\in [0,T]$, where $c_0=c_0(n,J[u_0,u_1,v_0,v_1])$ and $c^u_1,c^v_1$ are two positive constants independent of $T$.

To end this part, we recall the following propositions which are useful to prove Theorem \ref{Thm_Lower_Bound} in the next subsection.
\begin{prop}[Lemma $1$ in \cite{Matsumura}] \label{Estimates_Linear_Damped_Waves}
Let $n\geqslant 1$ and $k=0,1$. Then, the mild solutions to the linear Cauchy problem to \eqref{Eq_Linear_Damped_Waves} fulfill the following $(L^2(\mb{R}^n)\cap L^1(\mb{R}^n))-L^2(\mb{R}^n)$ and $L^2(\mb{R}^n)-L^2(\mb{R}^n)$ estimates:
\begin{align*}
\left\|\nabla^k w(t,\cdot)\right\|_{L^2(\mb{R}^n)} &\lesssim (1+t)^{-\frac{n}{4}- \frac{k}{2}}\left(\|w_0\|_{L^1(\mb{R}^n)}+ \|w_0\|_{H^k(\mb{R}^n)}+ \|w_1\|_{L^1(\mb{R}^n)}+ \|w_1\|_{L^2(\mb{R}^n)}\right), \\
\left\|\nabla^k w(t,\cdot)\right\|_{L^2(\mb{R}^n)} &\lesssim (1+t)^{- \frac{k}{2}}\left(\|w_0\|_{H^k(\mb{R}^n)}+ \|w_1\|_{L^2(\mb{R}^n)}\right).
\end{align*}
\end{prop}

\begin{prop}[The classical Gagliardo-Nirenberg inequality in \cite{ReissigEbert,Friedman}] \label{Classical_GN_In}
Let $r\in[1,\infty]$. It holds
\begin{align*}
	\|f\|_{L^r(\mb{R}^n)}\lesssim\|f\|_{L^2(\mb{R}^n)}^{1-\beta}\|\nabla f\|_{L^2(\mb{R}^n)}^{\beta}
\end{align*}
for $f\in \ml{C}_0^1(\mb{R}^n)$, where $\beta=n\left(\frac{1}{2}-\frac{1}{r}\right)$ and $\beta\in[0,1]$.
\end{prop}

\subsection{Lower bound estimates for the lifespan} \label{Sec.4.2}
In order to prove that the solution to \eqref{Eq_Coupled_Damped_Waves} satisfies the inequalities \eqref{Ineq.01} and \eqref{Ineq.02}, at first one deduces immediately the following estimate from the definition of the norm of $X(T)$ and Proposition \ref{Estimates_Linear_Damped_Waves}:
$$ \|(u^{\lin},v^{\lin})\|_{X(T)}\leqslant\varepsilon c_0(n,J[u_0,u_1,v_0,v_1]),$$
where is guaranteed by the fact
\begin{align*}
(1+t)^{-\gamma(p,q)}(\log(\mathrm{e}+t))^{\alpha(p,q)}\lesssim 1
\end{align*}
for any $t\in [0,T]$ due to $\gamma(p,q)>0$. Clearly, to achieve our aim, it suffices to demonstrate only
\begin{align}
\|u^{\non}\|_{Y_1(T)} &\leqslant c^u_1 (\log(\mathrm{e}+t))^{\alpha(p,q) }\,\|v\|_{Y_2(T)}^{p}, \label{Ineq.03} \\
\|v^{\non}\|_{Y_2(T)} &\leqslant c^v_1 (\log(\mathrm{e}+t))^{1-\alpha(p,q) q}\,\|u\|_{Y_1(T)}^{q}, \label{Ineq.04}
\end{align}
instead of \eqref{Ineq.01} and \eqref{Ineq.02}. Actually, by using the classical Gagliardo-Nirenberg inequality from Proposition \ref{Classical_GN_In}, we may derive 
\begin{align*}
	\|\,|u(\tau,\cdot)|^q\|_{L^1(\mb{R}^n)} &\lesssim (1+\tau)^{-\frac{n}{2}(q-1)+\gamma(p,q)q} (\log(\mathrm{e}+\tau))^{-\alpha(p,q) q}\|u\|_{Y_1(T)}^q \\
	&= (1+\tau)^{-1} (\log(\mathrm{e}+\tau))^{-\alpha(p,q) q}\|u\|_{Y_1(T)}^q,\\
	\|\,|u(\tau,\cdot)|^q\|_{L^2(\mb{R}^n)} &\lesssim (1+\tau)^{-\frac{n}{4}(2q-1)+\gamma(p,q)q} (\log(\mathrm{e}+\tau))^{-\alpha(p,q) q}\|u\|_{Y_1(T)}^q \\
	&= (1+\tau)^{-1-\frac{n}{4}} (\log(\mathrm{e}+\tau))^{-\alpha(p,q) q}\|u\|_{Y_1(T)}^q,
\end{align*}
and
\begin{align*}
	\|\,|v(\tau,\cdot)|^p\|_{L^1(\mb{R}^n)} &\lesssim (1+\tau)^{-\frac{n}{2}(p-1)}\|v\|_{Y_2(T)}^p= (1+\tau)^{-1+\gamma(p,q)}\|v\|_{Y_2(T)}^p,\\
	\|\,|v(\tau,\cdot)|^p\|_{L^2(\mb{R}^n)} &\lesssim (1+\tau)^{-\frac{n}{4}(2p-1)}\|v\|_{Y_2(T)}^p= (1+\tau)^{-1-\frac{n}{4}+\gamma(p,q)}\|v\|_{Y_2(T)}^p,
\end{align*}
for any $\tau\in[0,T]$. Let us sketch the proof of verification for these above estimates. On the one hand, we have utilized the relation from the critical curve
$$ \frac{q+1}{pq-1}=\frac{n}{2}\ \ \Rightarrow\ \ -\frac{n}{2}(q-1)+\gamma(p,q)q=-1,\ \ -\frac{n}{2}(p-1)=-1+\gamma(p,q) $$
in the powers of $(1+\tau)$. 
On the other hand, the following conditions must be satisfied due to the application of the classical Gagliardo-Nirenberg inequality:
\begin{align*}
 2\leqslant p,q\leqslant \ity  \ \ \mbox{if}\ \ n=1,2.
\end{align*}
 The first step is concerned with controlling the nonlinear integral terms $u^{\non}(t,\cdot)$ and $v^{\non}(t,\cdot)$ in the $L^2(\mb{R}^n)$ norm. By using the $(L^2(\mb{R}^n)\cap L^1(\mb{R}^n))-L^2(\mb{R}^n)$ estimate in $[0,t/2]$ and the $L^2(\mb{R}^n)-L^2(\mb{R}^n)$ estimate in $[t/2,t]$ from Proposition \ref{Estimates_Linear_Damped_Waves}, we obtain
\begin{align*}
	\|u^{\non}(t,\cdot)\|_{L^2(\mb{R}^n)} &\lesssim \int_0^{t/2}(1+t-\tau)^{-\frac{n}{4}}\|\,|v(\tau,\cdot)|^p\|_{L^2(\mb{R}^n)\cap L^1(\mb{R}^n)}\mathrm{d}\tau+ \int_{t/2}^t \|\,|v(\tau,\cdot)|^p\|_{L^2(\mb{R}^n)}\mathrm{d}\tau\\
	&\lesssim (1+t)^{-\frac{n}{4}}\,\|v\|_{Y_2(T)}^p\int_0^{t/2}(1+\tau)^{-1+\gamma(p,q)}\mathrm{d}\tau + (1+t)^{-1-\frac{n}{4}+\gamma(p,q)}\,\|v\|_{Y_2(T)}^p\int_{t/2}^t \mathrm{d}\tau\\
	&\lesssim (1+t)^{-\frac{n}{4}+\gamma(p,q)}\,\|v\|_{Y_2(T)}^p,
\end{align*}
where we employed the asymptotic relations $1+t-\tau \approx 1+t$ if $[0,t/2]$ and $1+\tau \approx 1+t$ if $[t/2,t]$ in the second inequality of the previous chain above. In other words, one gets
\begin{align}
	(1+t)^{-\gamma(p,q)+\frac{n}{4}} (\log(\mathrm{e}+t))^{\alpha(p,q) }\|u^{\non}(t,\cdot)\|_{L^2(\mb{R}^n)} &\lesssim (\log(\mathrm{e}+t))^{\alpha(p,q) }\,\|v\|_{Y_2(T)}^p. \label{Estimate_1}
\end{align}
The similar strategy to the previous one leads to
\begin{align*}
	\|v^{\non}(t,\cdot)\|_{L^2(\mb{R}^n)} &\lesssim \int_0^{t/2}(1+t-\tau)^{-\frac{n}{4}}\|\,|u(\tau,\cdot)|^q\|_{L^2(\mb{R}^n)\cap L^1(\mb{R}^n)}\mathrm{d}\tau + \int_{t/2}^t \|\,|u(\tau,\cdot)|^q\|_{L^2(\mb{R}^n)}\mathrm{d}\tau\\
	&\lesssim (1+t)^{-\frac{n}{4}}\,\|u\|_{Y_1(T)}^q\int_0^{t/2}(1+\tau)^{-1} (\log(\mathrm{e}+\tau))^{-\alpha(p,q) q}\mathrm{d}\tau\\
	&\quad+ (1+t)^{-1-\frac{n}{4}} (\log(\mathrm{e}+t))^{-\alpha(p,q) q}\,\|u\|_{Y_1(T)}^q\int_{t/2}^t \mathrm{d}\tau\\
	&\lesssim (1+t)^{-\frac{n}{4}} (\log(\mathrm{e}+t))^{1-\alpha(p,q) q}\,\|u\|_{Y_1(T)}^q
\end{align*}
due to the choice  the parameter $\alpha(p,q)$ satisfying
\begin{align*}
\alpha(p,q) q=\frac{n}{2}(p-1)=\frac{(q+1)(p-1)}{pq-1}=1-\gamma(p,q)<1.
\end{align*}
As a consequence, the following estimate holds:
\begin{align}
	(1+t)^{\frac{n}{4}}\|v^{\non}(t,\cdot)\|_{L^2(\mb{R}^n)}&\lesssim (\log(\mathrm{e}+t))^{1-\alpha(p,q) q}\,\|u\|_{Y_1(T)}^q. \label{Estimate_2}
\end{align}

In the second step, let us turn to the estimates for the gradient of solution. By repeating the same manner and analogous arguments as we dealt with $\|u^{\non}(t,\cdot)\|_{L^2(\mb{R}^n)}$ and $\|v^{\non}(t,\cdot)\|_{L^2(\mb{R}^n)}$, one also notices
\begin{align}
	(1+t)^{-\gamma(p,q)+\frac{n}{4}+\frac{1}{2}} (\log(\mathrm{e}+t))^{\alpha(p,q) }\|\nabla u^{\non}(t,\cdot)\|_{L^2(\mb{R}^n)} &\lesssim (\log(\mathrm{e}+t))^{\alpha(p,q) }\,\|v\|_{Y_2(T)}^p, \label{Estimate_3}\\
	(1+t)^{\frac{n}{4}+\frac{1}{2}}\|\nabla v^{\non}(t,\cdot)\|_{L^2(\mb{R}^n)}&\lesssim (\log(\mathrm{e}+t))^{1-\alpha(p,q) q}\,\|u\|_{Y_1(T)}^q. \label{Estimate_4}
\end{align}
All in all, by the definition of the corresponding norms in $Y_1(T)$ and $Y_2(T)$ we can link the obtained estimates from \eqref{Estimate_1} to \eqref{Estimate_4} to claim \eqref{Ineq.03} and \eqref{Ineq.04} automatically. 

Afterwards, motivated by the approach in \cite{Iked-Ogaw-2016} we determine
\begin{align*}
	T^*:=\sup\left\{T\in[0,T_{\varepsilon,\mathrm{m}})\,\, \text{ such that }\,\, F(T):= \|(u,v)\|_{X(T)}\leqslant M\varepsilon \right\}
\end{align*}
with a sufficiently large constant $M>0$, which will be defined in next steps. Thanks to the fact $\|u\|_{Y_1(T^*)}\leqslant \|(u,v)\|_{X(T^*)}\leqslant M\varepsilon$, it holds from \eqref{Ineq.02} that
\begin{align}
	\|v\|_{Y_2(T^*)}\leqslant   \varepsilon c_0+ c^v_1 M^q (\log(\mathrm{e}+T^*))^{1-\alpha(p,q) q}\varepsilon^q, \label{Estimate_5}
\end{align}
where the restriction $1-\alpha(p,q) q>0$ was used again. Then, substituting the estimate \eqref{Estimate_5} into the inequality \eqref{Ineq.01} results
\begin{align*}
\|u\|_{Y_1(T^*)}&\leqslant \varepsilon c_0+ (\log(\mathrm{e}+T^*))^{\alpha(p,q) }\left(c_2 \varepsilon^{p}+ c_3 M^{pq}(\log(\mathrm{e}+T^*))^{p(1-\alpha(p,q) q)}\varepsilon^{pq}\right)\\
&\leqslant \varepsilon\left(c_0+ c_2(\log(\mathrm{e}+T^*))^{\alpha(p,q)}\varepsilon^{p-1}+ c_3 M^{pq}(\log(\mathrm{e}+T^*))^{p-\alpha(p,q) (pq-1)}\varepsilon^{pq-1} \right)
\end{align*}
with two positive constants $c_2= c_2(c_0,c^u_1,p)$ and $c_3= c_3(c^u_1,c^v_1,p)$. Here, we recall $\alpha(p,q)>0$ to deduce the increasing property of the logarithmic function. So, we can take a large constant $M>0$ such that $0<c_0 < M/8$ to enjoy
\begin{align*}
	\|u\|_{Y_1(T^*)}&< \frac{3}{8}M\varepsilon,
\end{align*}
providing that
\begin{align*}
	8c_2 M^{-1}(\log(\mathrm{e}+T^*))^{\alpha(p,q)}\varepsilon^{p-1}< 1\ \ \mbox{and}\ \ 8c_3 M^{pq-1}(\log(\mathrm{e}+T^*))^{p-\alpha(p,q) (pq-1)}\varepsilon^{pq-1}< 1.
\end{align*}
What's more, it holds from \eqref{Estimate_5} that
\begin{align*}
	\|v\|_{Y_2(T^*)}\leqslant \varepsilon\left(c_0+ c^v_1 M^q (\log(\mathrm{e}+T^*))^{1-\alpha(p,q) q}\varepsilon^{q-1}\right)< \frac{1}{4}M\varepsilon,
\end{align*}
when we consider
\begin{align*}
	8c^v_1 M^{q-1} (\log(\mathrm{e}+T^*))^{1-\alpha(p,q) q}\varepsilon^{q-1}< 1.
\end{align*}
Collecting the above two estimates, we know
\begin{align}\label{Ineq_Contri}
	F(T^*)= \|(u,v)\|_{X(T^*)}= \|u\|_{Y_1(T^*)}+\|v\|_{Y_2(T^*)}< \frac{5}{8}M\varepsilon< M\varepsilon.
\end{align}
Note that $F=F(T)$ is a continuous function for any $T\in(0,T_{\varepsilon,\mathrm{m}})$. Thus, it follows from \eqref{Ineq_Contri} that there exists a time $T_0\in(T^*,T_{\varepsilon,\mathrm{m}})$ satisfying $F(T_0)\leqslant M\varepsilon$ so that the contradiction appears to the definition of $T^*$. In other words, we have to pose that one of the following estimates are true:
\begin{align*}
	&8c_2 M^{-1}(\log(\mathrm{e}+T^*))^{\alpha(p,q)}\varepsilon^{p-1}\geqslant1,\\
	&8c_3 M^{pq-1}(\log(\mathrm{e}+T^*))^{p-\alpha(p,q) (pq-1)}\varepsilon^{pq-1}\geqslant 1,\\
	&8c^v_1 M^{q-1}(\log(\mathrm{e}+T^*))^{1-\alpha(p,q) q}\varepsilon^{q-1}\geqslant 1.
\end{align*}
Then, we can find the blow-up time
\begin{align*}
	T_{\varepsilon,\mathrm{m}}\geqslant\exp\left(c \varepsilon^{-\min\left\{ \frac{p-1}{\alpha(p,q)},\frac{pq-1}{p-\alpha(p,q)(pq-1)},\frac{q-1}{1-\alpha(p,q) q}  \right\}}\right)=
		\exp\big(c \varepsilon^{- \frac{p-1}{\alpha(p,q)}}\big)
\end{align*}
 where $c$ is a positive constant independent of the small parameter $\varepsilon$. Here, we pay attention to the condition $p-\alpha(p,q) (pq-1)>0$ when $\alpha(p,q) q<1$ holds. Thus, we may rewrite the lower bound estimates by
 \begin{align*}
 	T_{\varepsilon,\mathrm{m}}\geqslant \mathrm{exp}\left(c\e^{-\frac{q(pq-1)}{q+1}}\right).
 \end{align*}
 By this way, we can achieve our aim to show the really sharp lifespan in the case $p<q$. 
 
%

The special case $p=q$ would be dealt more simply (without any logarithmic weighted function) with an analogous procedure to the proof of the case $p<q$ by setting $\gamma(p,q)=0$ and $\alpha(p,q)=0$. More precisely, from \eqref{Estimate_5} one arrives at
\begin{align*}
	\|v\|_{Y_2(T^*)}\leqslant \varepsilon \left(c_0+ c^v_1 M^p \log(\mathrm{e}+T^*) \varepsilon^{p-1}\right).
\end{align*}
Again, following some arguments used in \cite{Iked-Ogaw-2016} we obtain
\begin{align*}
	T_{\varepsilon,\mathrm{m}}\geqslant \exp\left(c\varepsilon^{-(p-1)}\right).
\end{align*}
 Hence, our proof is completed.

\section{Concluding remarks}\label{Section_Concluding_Remarks}
\begin{remark}
	\fontshape{n}
	\selectfont
We will give an application of our methods to the weakly coupled system of semilinear reaction-diffusion equations in this remark. As we know, this system can describe a model for heat propagations in a two-component combustible mixture. More precisely, let us consider the following Cauchy problem for two scalar functions $\tilde{u}=\tilde{u}(t,x)$, $\tilde{v}=\tilde{v}(t,x)$ standing for the temperatures of the interacting components:
\begin{align}\label{Eq_Coupled_Heats}
\begin{cases}
\tilde{u}_{t}-\Delta \tilde{u}=|\tilde{v}|^p,&x\in\mb{R}^n,\ t\in(0,T),\\
\tilde{v}_{t}-\Delta \tilde{v}=|\tilde{u}|^q,&x\in\mb{R}^n,\ t\in(0,T),\\
(\tilde{u},\tilde{v})(0,x)=(\varepsilon \tilde{u}_0,\varepsilon \tilde{v}_0)(x),&x\in\mb{R}^n,
\end{cases}
\end{align}
where $p,q>1$ satisfy the critical condition $\alpha_{\max}(p,q)=n/2$ for any $n\geqslant 1$ (see \cite{Esco-Herr-1991,Esco-Levi-1995,Moch-Huan-1998,Ishi-Kawa-Sie-1998,Renc-2000,Umed-2003,Aoya-Tsut-Yama-2007} and references therein, particularly, a complete introduction in \cite{Fuji-Iked-Waka-2020}). Here, thermal conductivity is supposed constant and equal for both substances. For one thing, by following the same approach as the proof of Theorem \ref{Thm_Upper_Bound}, we also can conclude the same upper bound estimates for the lifespan of solutions to \eqref{Eq_Coupled_Heats}. We assume that $\tilde{u}_0,\tilde{v}_0\in\ml{C}_0^{\infty}(\mb{R}^n)$ enjoy
$$\int_{\mb{R}^n}\tilde{u}_0(x)\mathrm{d}x>0 \ \  \text{as well as}\ \  \int_{\mb{R}^n}\tilde{v}_0(x)\mathrm{d}x>0. $$
Then, the upper bound estimates for the lifespan $\widetilde{T}_{\varepsilon}$ of weak solutions to \eqref{Eq_Coupled_Heats} under the situation $\alpha_{\max}(p,q)=n/2$ fulfill the same estimates as \eqref{Lifespan}. For another thing, since our approach to deal with the lower bound estimates depends on some decay estimates for solutions, we expect that following the same procedure as the proof of Theorem \ref{Thm_Lower_Bound} associated with some well-known $L^{r}(\mb{R}^n)-L^{m}(\mb{R}^n)$ estimates with $1\leqslant m\leqslant r\leqslant \infty$ allows one to derive the same lower bound estimates for the lifespan as \eqref{Lifespan_Lower} to \eqref{Eq_Coupled_Heats} in the critical case $\alpha_{\max}(p,q)=n/2$. Together with these expectations, we may believe that the sharp lifespan estimates $\widetilde{T}_{\varepsilon}$ for solutions to the reaction-diffusion systems \eqref{Eq_Coupled_Heats} in low spatial dimensions with $\alpha_{\max}(p,q)=n/2$ are defined by
\begin{align*}
\widetilde{T}_{\varepsilon}\sim\begin{cases}
\exp\left( C\varepsilon^{-(p-1)} \right)&\mbox{if}\ \ p=q,\\
\exp\left(C\e^{-(pq-p_{\mathrm{Fuj}}(n))}\right) &\mbox{if}\ \ p\neq q,
\end{cases}
\end{align*}
where $C>0$ is a constant independent of $\varepsilon$.
\end{remark}
\begin{remark}
	\fontshape{n}
	\selectfont
Throughout this paper, we have succeeded in deriving not only some upper bound estimates but also some lower bound estimates of solutions to the Cauchy problem for the weakly coupled system of semlinear damped wave equations in the critical case. Clearly, this is to state that the obtained lifespan estimates in this work are actually sharp. More generally, we expect that our approach utilized in this paper can be applied to study lifespan estimates of solutions to the Cauchy problem for other weakly coupled systems of semilinear parabolic-like evolution equations in the critical case, for example, weakly coupled systems of semilinear wave equations with time-dependent damping terms \cite{Djau-Reis-2018,Djau-Reis-2019,Chen-Palmieri-2019}, or weakly coupled systems of semilinear $\sigma$-evolution equations with damping terms \cite{D'Abbicco2015,Dao-2019} for some suitable parameters.
\end{remark}

\section*{Acknowledgments}
 This research of the second author (Tuan Anh Dao) is funded (or partially funded) by the Simons Foundation Grant Targeted for Institute of Mathematics, Vietnam Academy of Science and Technology. The authors thank Alessandro Palmieri (University of Pisa) and Ya-guang Wang (Shanghai Jiao Tong University) for their suggestions in the preparation of the paper.

\end{document}